\documentclass[a5paper,twoside]{article}
\usepackage[utf8]{inputenc}
\usepackage{lmodern}
\usepackage{amsmath}
\usepackage{amssymb}
\usepackage{cmap}
\usepackage[T2A,T1]{fontenc}
\usepackage[russian,english,italian,german,french]{babel}
\newcommand\italian[1]{\foreignlanguage{italian}{#1}}
\newcommand\german[1]{\foreignlanguage{german}{#1}}
\newcommand\english[1]{\foreignlanguage{english}{#1}}
\usepackage[outer=10mm,inner=12mm,top=12mm,bottom=15mm,footskip=7mm]{geometry}
\usepackage{paralist}
\makeatletter
\renewcommand\section{\@startsection{section}{1}{\parindent}{\medskipamount}{\medskipamount}{\normalfont\normalsize\fontseries{b}\selectfont}}
\renewcommand\subsection{\@startsection{subsection}{2}{\parindent}{\medskipamount}{\medskipamount}{\normalfont\normalsize\fontshape{it}\selectfont}}
\makeatother

\usepackage{titlesec}
\titlelabel{\thetitle.\ }
\usepackage{doi}
\PassOptionsToPackage{hyphens}{url}
\usepackage{hyperref}
\hypersetup{colorlinks=true, linkcolor=blue, citecolor=blue, filecolor=blue, urlcolor=blue}

\usepackage{breakurl}

\title{\vspace{-.5in}Roger Apéry, l'humour au service d'une pensée libre et originale sur les mathématiques constructives.}
\author{Henri Lombardi \and Stefan Neuwirth}
\date{}

\begin{document}

\maketitle

\begin{sloppypar}
Cet article invite à la lecture d'un texte célèbre de Roger Apéry, «Mathématique constructive», paru dans l'ouvrage \emph{Penser les mathématiques: séminaire de philosophie et mathématiques de l'École normale supérieure (J.~Dieudonné, M.~Loi, R.~Thom)} édité par F.~Guénard et G.~Lelièvre, Paris: Seuil, 1982, p.~58-72. Il est disponible sous \url{https://hal.archives-ouvertes.fr/hal-01522168}.
\end{sloppypar}

Le texte original de la conférence donnée en 1976 au «séminaire Loi» est disponible sur Numdam sous \url{http://www.numdam.org/item?id=SPHM_1976___1_A1_0}.

François Apéry a écrit une courte note biographique de son père pour la revue \emph{\english{The Mathematical Intelligencer}} en 1996. Une traduction française est disponible sous \url{http://peccatte.karefil.com/PhiMathsTextes/AperyFR.htm}.

\section*{Notre lecture.}

En l'an 1976, en pleine période de bourbakisme triomphant, Roger Apéry est invité à donner une conférence au séminaire de philosophie et mathématiques de l'École normale supérieure. Jean Dieudonné, représentant officiel du bourbakisme parmi les organisateurs du Séminaire, est un ami personnel de Roger Apéry, comme le rappelle la note biographique citée. Ainsi, malgré le slogan «À bas Euclide» (appuyé par un célèbre pensum de géométrie sans figures ni intuition publié par Dieudonné), la guerre n'est pas totale et les adeptes du paradis de Cantor dans la version formaliste tolèrent que l'on développe une philosophie des mathématiques en s'opposant au cadre officiel dominant. Et pourtant Bourbaki avait déclaré la guerre terminée par triomphe complet et autoproclamé des vainqueur\textperiodcentered e\textperiodcentered s.

Le texte est empli d'un humour parfois dévastateur, comme la fin du premier paragraphe.
\begin{quotation}
  \noindent Enfin, aucun argument solide ne permet d'affirmer que L.~Kronecker, H.~Poincaré ou H.~Weyl étaient plus angoissés que Cantor, Hilbert ou Russell.
\end{quotation}

Le texte de la conférence est parfois plus direct et mordant que le texte remanié paru dans \emph{Penser les mathématiques}, même s'il n'y a aucune variation dans la pensée philosophique.

À titre d'exemple, voici le passage décrivant la position du formalisme dans le texte de la conférence. Beaucoup moins argumenté dans les détails que le texte remanié, beaucoup plus court, il n'en est que plus incisif.
\begin{quotation}
Soucieux d'éliminer les antinomies issues du cantorisme et d'extirper toute trace de subjectivisme, les formalistes réduisent la mathématique à un simple jeu défini par des règles. Ils abandonnent dédaigneusement au psychologue l'activité mentale du mathématicien dont l'examen leur semble inutile pour juger la validité des résultats. L'adéquation au monde physique leur parait due à une «harmonie préétablie» ou à un «miracle».  Ils ne veulent connaitre des mathématiques, de la musique, de la littérature que leur formalisation écrite conservée par les bibliothèques.

Attachés sentimentalement au «paradis créé pour nous par Cantor», ils pratiquent le double langage: quand le profane croit qu'ils démontrent des vérités, ils se contentent d'établir la possibilité d'obtenir leurs résultats en respectant les règles d'un jeu formel, de façon analogue au problème d'échecs énonçant: les blancs jouent et gagnent.

Leur hostilité à la liberté du sujet les pousse à la création d'un dieu mathématique à plusieurs personnes, qui tente d'être immortel en renouvelant périodiquement ses membres; ce dieu \emph{révèle} aux populations les bonnes définitions et les bonnes théories.

En face de ceux qui refusent la science au nom du bon sens, de la poésie, de la religion, les formalistes rejettent comme démunis de sens les concepts d'espace, de temps, de liberté.
\end{quotation}

Le dieu mathématique à plusieurs personnes qui tente d'être immortel est bien évidemment le groupe Bourbaki, et il occupe dans le discours d'Apéry la position symétrique des gens qui refusent la science au nom de la religion.  Il écrit ceci après avoir montré que Leibniz et Cantor avaient besoin de Dieu pour se convaincre de croire en l'infini actuel. Voici donc les platonicien\textperiodcentered ne\textperiodcentered s et les formalistes renvoyé\textperiodcentered e\textperiodcentered s au paradis, bien loin de la science.  Naturellement, la plupart des mathématicien\textperiodcentered ne\textperiodcentered s sont plutôt platonicien\textperiodcentered ne\textperiodcentered s que formalistes, car il\textperiodcentered elle\textperiodcentered s voient bien que les mathématiques ne se réduisent ni à un plan de carrière de mathématicien\textperiodcentered ne ni à un jeu formel qui, par «miracle», s'accorde avec le monde réel objectif\footnote{Comme le montrent la mécanique de Newton, les équations de Maxwell, les lentilles de Fresnel, la prévision météo, le vol des avions, le lancement des satellites, les technologies modernes, les disques compacts, les téléphones portables ou non, le GPS, et, hélas!, les bombes à fragmentation, la bombe atomique («et c'est depuis qu'ils sont civilisés\dots»), les drones et les mathématiques financières.}.

L'ironie mordante d'Apéry se manifeste aussi contre la réforme des mathématiques modernes, qui s'avèrera effectivement catastrophique. C'est le point~10\up{o} du formalisme dans le texte remanié.
\begin{quotation}
Uniformiser les esprits par l'enseignement des «mathématiques modernes», où on laisse croire aux enfants qu'entourer des petits objets par une ficelle est une activité mathématique au lieu de leur apprendre à compter, à calculer et à examiner les propriétés des figures.
\end{quotation}

L'aspect le plus stimulant du texte d'Apéry est la comparaison étroite qu'il établit entre les mathématiques et d'autres activités humaines comme la littérature et la musique, qui n'existent qu'avec le déroulement du temps et eu égard à la psychologie humaine. En insistant sur le fait que «celui qui possède des textes mathématiques dont il ne comprend pas l'articulation ne possède rien», il nous rappelle une évidence profonde que nous oublions quand nous prenons les mathématiques comme une sorte de vérité absolue préétablie révélée et non soumise à discussion.

Que sont donc les objets mathématiques? Les platonicien\textperiodcentered ne\textperiodcentered s les conçoivent comme préexistants de toute éternité dans un monde idéal des idées, la plupart du temps ensembliste de surcroit.
\begin{quotation}
Comme le platonicien et contrairement au formaliste, le mathématicien constructif reconnait une certaine réalité aux objets mathématiques, mais les différencie essentiellement des objets matériels, en ne leur attribuant que les propriétés susceptibles de démonstration. Une distinction analogue différencie les héros de roman des personnages historiques.
\end{quotation}

Suit une comparaison savoureuse entre les statuts de vérité différents concernant Vercingétorix et Don Quichotte, aboutissant à cette magnifique chute:
\begin{quotation}
  \noindent l'ensemble des réels, comme Don Quichotte, est un être essentiellement incomplet.
\end{quotation}

Sous une forme particulièrement condensée, la conception des objets mathématiques d'Apéry est résumée comme suit.
\begin{quotation}
Selon la conception constructive, il n'y a pas de mathématique sans mathématicien. En tant qu'êtres de raison, les êtres mathématiques n'existent que dans la pensée du mathématicien et non dans un monde platonicien indépendant de l'esprit humain [\dots].
\end{quotation}

Cela ne conduit cependant pas Apéry jusqu'à la position intuitionniste de Brouwer. Voici comment il s'en explique dans le texte initial de la conférence.
\begin{quotation}
En face des «statiques» qui veulent détruire l'intuition (on connait les résultats désastreux dans l'enseignement), Brouwer lui laisse une trop large part en considérant comme prouvé ce qui est intuitivement clair; la clarté intuitive varie souvent d'un mathématicien à l'autre.

\begin{sloppypar}
Nous préférons au vocable «intuitionnisme» utilisé par Brouwer, le vocable «constructivisme» qui évoque mieux les méthodes de preuve permises.
\end{sloppypar}

\end{quotation}

En résumé, Apéry se situe dans les courants des mathématiques constructives représentés à son époque par Bishop et Shanin (cités à la fin de son texte), sans prendre parti entre leurs différences d'approche.

La section finale dans laquelle il explique très clairement pourquoi la logique constructive est plus riche que la logique classique est une œuvre de salubrité publique qui a dû résonner de manière bizarre aux oreilles de Dieudonné. Cette réalité objective est toujours mal assimilée (et souvent simplement ignorée) par la plupart des mathématicien\textperiodcentered ne\textperiodcentered s. Les théoricien\textperiodcentered ne\textperiodcentered s de l'informatique sont en général mieux informé\textperiodcentered e\textperiodcentered s.

L'idéal constructif poursuivi par Gauss, Kronecker et Bishop, visant à expliciter tout théorème de mathématiques sérieux comme un algorithme, est devenu beaucoup plus crédible et concret avec l'avènement des calculs sur machine, inaccessibles aux humains mais con\-trô\-lés par eux. Presque tou\textperiodcentered te\textperiodcentered s les mathématicien\textperiodcentered ne\textperiodcentered s considèrent aujourd'hui qu'une version explicite d'un théorème présente un réel intérêt par rapport à une version purement idéale.

Peut-être est-ce ici un bon endroit pour expliquer ce qui distingue les mathématiques constructives «minimales» à la Bishop des mathématiques constructives à la Shanin. En très bref et avec certainement un peu d'exagération, les mathématiques de Bishop sont un cadre commun de base pour toutes les mathématiques pratiquées aujourd'hui. Les objets mathématiques sont essentiellement les mêmes que pour les mathématiques qui ajoutent des axiomes sujets à discussion.

Les mathématiques classiques ajoutent le principe du tiers exclu et l'axiome du choix. Le principe du tiers exclu est le ver dans le fruit qui empêche presque toujours que les démonstrations classiques recèlent un contenu algorithmique évident. Un travail de décryptage est presque toujours nécessaire, parfois facile, parfois fort délicat\footnote{Les auteurs de cet article estiment que ce travail n'est pas reconnu à sa juste valeur par les structures de la recherche et peuvent témoigner qu'il n'a aucune incidence positive sur la carrière (tout au contraire il la ralentit). Cela explique pourquoi il n'est que rarement entrepris par les mathématiciens.}.

Les mathématiques à la Shanin affirment une «thèse de Church-Turing renforcée» selon laquelle tout objet mathématique est codable en machine de Turing et tous les théorèmes de mathématiques relèvent de programmes de calculs à la Turing. Cette position radicale conduit à des ennuis sérieux avec des théorèmes peu intuitifs d'une part et en contradiction immédiate avec les mathématiques classiques d'autre part. Par exemple toute fonction de $\mathbb{R}$ dans $\mathbb{R}$ est continue en tout point, mais il y a des fonctions uniformément continues sur l'intervalle [0,1] partout strictement positives avec une borne inférieure nulle. C'est une rançon cher payée pour le fait de n'admettre à priori et à tout jamais que les réels «calculables sur machine». Mieux vaut penser que les réels ressemblent plus à Don Quichotte qu'à Vercingétorix.

En fait la source de divergence la plus importante entre Bishop et Shanin porte sur la question de la notion même de construction. Cela se manifeste dès les suites calculables de nombres entiers. Shanin pense que ce concept est bien défini dans l'absolu, par le recours aux machines de Turing. Bishop dit que l'on a là affaire à une notion de base, la notion de construction, non susceptible d'être définie et sur laquelle les mathématicien\textperiodcentered ne\textperiodcentered s doivent s'accorder au cas par cas. Quand un\textperiodcentered e mathématicien\textperiodcentered ne imagine une suite calculable d'entiers, ce n'est généralement pas en termes de programme de Turing. Même si en pratique il\textperiodcentered elle arrive en fin de compte à en déduire un programme, il s'agit d'un nouveau travail qui n'a aucun caractère mécanique. Les calculs mécaniques à la Turing ne sont donc pas à priori la même chose que les constructions au sens intuitif et indéfini de la chose. En outre, et c'est là un point très important, affirmer qu'un programme donné produira bien un calcul fini qui terminera un jour demande une démonstration constructive qui ne peut en aucun cas être considérée comme mécanisable. Cela résulte du fameux «théorème de l'arrêt» de Turing, qui n'est autre que la version machine de Turing de la méthode diagonale de Cantor.

\section*{Une bibliographie commentée.}

Nous donnons ci-dessous une liste de références en rapport avec notre propos.
Le\textperiodcentered la lecteur\textperiodcentered rice trouvera aussi quelques références purement techniques à la fin de la note technique qui suit.
Signalons qu'une bibliographie commentée est déjà parue à la fin de l'article «Le programme de Hilbert et les mathématiques constructives», \emph{Repères IREM} \no50, 2003, p.~85-104, \url{http://www.univ-irem.fr/spip.php?rubrique24&id_numero=50}.

\subsection*{Voici tout d'abord deux textes qui parlent de Roger Apéry.}

{\leftskip15pt
\parindent-15pt

P.~Ageron, «La philosophie mathématique de Roger Apéry», \emph{Philosophia Scientiæ}, cahier spécial~5, 2005, p.~233-256, \url{http://philosophiascientiae.revues.org/89}.

}L'auteur résume ainsi son article:
\begin{quotation}
\noindent Pour qui s'intéresse à la philosophie des mathématiques, Roger Apéry (1916-1994) incarne le défenseur de la mathématique constructive et l'adversaire résolu du formalisme et du bourbakisme. On sait moins qu'il est aussi l'un des premiers universitaires français à avoir fait la promotion de la théorie des catégories, pourtant hautement structuraliste et souvent jugée comme très formelle. L'objectif principal de notre étude est de préciser les conditions historiques et la teneur philosophique du double enthousiasme d'Apéry, afin de vérifier la cohérence d'une pensée libre, originale et attachante.
\end{quotation}

{\leftskip15pt
\parindent-15pt
\sloppy

F.~Apéry, «Roger Apéry, \english{1916-1994: a radical mathematician}», \emph{\english{The Mathematical Intelligencer}}, vol.~18, \no2, 1996, p.~54-61, \url{http://peccatte.karefil.com/PhiMathsTextes/Apery.html}, traduction par P.~Karila et M.~Saunier: \url{http://peccatte.karefil.com/PhiMathsTextes/AperyFR.htm}.

}\subsection*{Ensuite des textes de Roger Apéry lui-même.}

{\leftskip15pt
\parindent-15pt

R.~Apéry, «Axiomes et postulats», dans \emph{Actes du X\up{me} congrès international de philosophie (Amsterdam, 11-18 aout 1948)}, vol.~II, Amsterdam: North-Holland, 1949, p.~708-710.

R.~Apéry, «Le rôle de l'intuition en mathématiques», dans \emph{Congrès international de philosophie des sciences (Paris, 17-22 octobre 1949)}, vol.~III, Paris: Hermann \& C\up{ie}, 1951, p.~85-88.

------, «Les mathématiques sont-elles une théorie pure?», \emph{Dialectica}, vol.~6, 1952, p.~309-310.

}\subsection*{\sloppy Enfin quelques textes importants pour les mathématiques constructives, avec quelques commentaires.}

Quand Errett Bishop publie en 1967 le livre \emph{\english{Foundations of Constructive Analysis}}, où il interprète en termes constructifs les bases de l'analyse moderne, il réalise un morceau substantiel du programme de Hilbert relu sous la forme suivante:
\begin{asparaitem}
\item lorsqu'un résultat concret est démontré en mathématiques par des méthodes douteuses, certifier ce résultat par des méthodes sûres;
\item réaliser ce travail de manière aussi systématique (voire automatique) que possible.
\end{asparaitem}
Dans le livre de Bishop, tous les théorèmes d'analyse ont la signification d'algorithmes qui calculent des objets concrets à partir d'autres objets concrets, conformément à certaines spécifications requises, et ces algorithmes sont prouvés par des méthodes sûres: en particulier personne ne conteste qu'ils aboutissent certainement en un temps fini au résultat souhaité. Ainsi les bases de l'analyse sont ramenées à un degré de certitude comparable à celui qui règne en théorie élémentaire des entiers naturels.

Bishop va bien au delà de ce qu'avait pu faire auparavant un logicien remarquable comme Goodstein: non seulement sont traités une quantité de résultats incomparablement plus grande, mais encore, le style d'exposition est direct, sans autre différence sensible avec le style mathématique usuel qu'une attention scrupuleuse accordée aux aspects effectifs.

On pourra lire à ce sujet l'article de D.~Knuth dans lequel il analyse la page~100 de différents livres de mathématiques, dont celui de Bishop, du point de vue de la pensée algorithmique.

Non seulement, le programme de Hilbert (revisité) n'est pas utopique, mais il a de bonnes chances de pouvoir être développé en grand après un tel coup de maitre.

Mais au lieu d'être acclamé comme le travail d'un bienfaiteur des mathématiques, ce livre a été accueilli par une quasi-indifférence des professionnels. Ceci ne s'explique que par le poids de la routine (qu'est-ce que c'est que ce type qui prétend faire changer nos standards?) et par le manque presque total de questionnement des professionnels concernant la signification de leur activité scientifique. L'épistémologie des mathématiques n'est pas à l'ordre du jour, elle ne fait pas partie du cursus normal: elle n'est presque pas enseignée, et quand elle l'est c'est en général uniquement à titre de spécialité. Le livre de Bishop a été rapidement épuisé. Il a été réimprimé récemment. Il a fait l'objet d'études approfondies chez les logicien\textperiodcentered ne\textperiodcentered s (voir par exemple les textes de M.~Beeson). Une deuxième édition en collaboration avec D.~Bridges, dans laquelle la théorie de l'intégrale de Lebesgue a été modifiée, est parue en 1985.

En algèbre, le point de vue algorithmique a toujours eu des défenseurs. Il y a de quoi, puisque le mot «algèbre» est tiré de «al-jabr», extrait du titre d'un livre écrit il y a fort longtemps par un auteur qui s'appelait «M.~Algorithme» (Al-Khwarizmi). Il faut bien évidemment souligner la tradition de Gauss et Kronecker, entièrement dans le style algorithmique. Bien que les méthodes abstraites soient ensuite devenues quelque peu hégémoniques sous l'influence de Hilbert puis de Bourbaki, il est encore fréquent d'enseigner et de publier des algorithmes. Signalons entre autres les travaux de Seidenberg. En 1988, le merveilleux petit livre de Mines, Richman et Ruitenburg a fait pour les bases de l'algèbre moderne ce qu'avait fait le livre de Bishop pour celles de l'analyse.

La nouvelle discipline du calcul formel (calculs symboliques et algébriques sur machine) se rattache de facto à cette tradition, même si les auteur\textperiodcentered e\textperiodcentered s ne comprennent pas toujours la différence entre preuve constructive et algorithme: un algorithme qui met en œuvre un théorème d'algèbre est parfois prouvé par des méthodes abstraites non constructives, auquel cas le programme de Hilbert n'est réalisé que très imparfaitement pour le théorème en question. Voir le livre édité par Cohen, Cuypers et Sterk et celui de Cox, Little et O'Shea.

Le livre de Bridges et Richman est une très bonne introduction aux différentes variantes des mathématiques constructives.

Le livre d'Aberth est une présentation simple de la variante récursiviste russe (dans laquelle tous les objets
sont supposés récursifs). 

Le livre d'Ageron est agréable à lire: le style est à la fois simple, rigoureux et informel. Il expose quelques fondements de la théorie des ensembles et de celle des catégories sans le principe du tiers exclu. De nombreuses explications historiques éclairent le propos. Une dose d'humour est aussi présente, comme dans les trois solutions proposées pour s'en sortir face à l'absence de l'ensemble de tous les ensembles: a)~jouer avec les mots; b)~proclamer: il est interdit d'interdire; c)~vivre avec ses contradictions.

Le livre de David, Nour et Raffalli est un livre d'enseignement universitaire. Il doit être salué comme le premier livre de ce type (à notre connaissance) écrit en français et qui accorde la place qu'elle mérite à la logique intuitionniste (la logique des mathématiques constructives). Cependant il est écrit du point de vue des mathématiques classiques, ce qui est paradoxal quand on veut traiter sérieusement la question des fondements (il est vrai que ce n'est pas l'objet essentiel de l'ouvrage).

Le livre de Gilles Dowek est un excellent petit texte de présentation de la logique, mais écrit pour un public large, donc sans description précise des outils de la logique mathématique.

Le livre de Jean-Louis Krivine, d'une clarté remarquable, doit être lu par tout\textperiodcentered e mathématicien\textperiodcentered ne qui veut comprendre ce qu'est la théorie formelle des ensembles: l'étude d'une structure particulière, non pas celle d'anneau commutatif, ni celle de treillis distributif, mais celle d'univers, c'est-à-dire la structure des ensembles naïfs munis d'une relation notée «$\in$» vérifiant un système d'axiomes mis au point pour éviter le paradoxe de Russell.

En informatique, les auteur\textperiodcentered e\textperiodcentered s clairvoyant\textperiodcentered e\textperiodcentered s ne cachent pas que la seule logique qui vaille pour l'informatique théorique est la logique intuitionniste. Voir par exemple le beau livre de R.~Lalement.

Le livre de Lakatos reste une source importante de réflexion sur les fondements. Notamment, la place centrale qui est accordée aux preuves par rapport aux théorèmes est en accord profond avec le dicton constructif: ce qui est vrai est ce qui peut être prouvé.

L'article d'Abraham Robinson, le fondateur de l'analyse non standard, montre à quel point ce visionnaire des infinitésimaux doute de la réalité de l'infini cantorien.

Le livre de Nagel, Newman, Gödel et Girard donne une traduction française de l'article original de Gödel avec le théorème d'incomplétude. Il contient aussi la traduction française d'un texte de Nagel et Newman exposant de manière simple les idées à l'œuvre dans l'article de Gödel. Enfin on trouve un commentaire de J.-Y.~Girard dans lequel il assassine avec brio un certain nombre de positions adverses, dont celle de Hilbert.

Le livre édité par Toraldo di Francia est un recueil d'articles en français, en italien et en anglais de physicien\textperiodcentered ne\textperiodcentered s, de logicien\textperiodcentered ne\textperiodcentered s, de philosophes et d'historien\textperiodcentered ne\textperiodcentered s des sciences sur le problème de l'infini dans les sciences.  Sa lecture pourrait aider les mathématicien\textperiodcentered ne\textperiodcentered s à sortir de leur bulle et à jeter un regard critique sur ce qui leur semble évident par habitude professionnelle.

Le livre de Feferman reprend une série d'articles importants concernant les fondements des mathématiques. Tout son livre tend à la conclusion que pour faire des mathématiques, la théorie des ensembles ne sert à rien. Elle sert évidemment à développer des intuitions fructueuses, mais, sur le fond, elle n'apporte aucun outil que ne nous apporterait pas la théorie des entiers naturels. Le point de vue personnel de Feferman semble être essentiellement un développement de la position défendue par Hermann Weyl dans \emph{Le continu}. Un des systèmes formels proposés par Feferman essaie de traduire au plus près les intuitions de Weyl, et ce système est «complètement sûr» puisque c'est une extension conservative de la théorie formelle de Peano pour les entiers naturels.\medskip

{\leftskip15pt
\parindent-15pt

O.~Aberth, \emph{\english{Computable analysis}}, New York: McGraw-Hill, 1980.

P.~Ageron, \emph{Logiques, ensembles, catégories: le point de vue constructif}, Paris: Ellipses, 2000.

M.~Beeson, \emph{\english{Foundations of constructive mathematics}}, Berlin: Springer, 1985.

------, «\english{Some theories conservative over intuitionistic arithmetic}», dans \emph{\english{Logic and computation: proceedings of a workshop held at Carnegie Mellon University, June 30-July 2, 1987}} (éd. par W.~Sieg), Providence: \english{American Mathematical Society}, 1990, p.~1-15.

E.~Bishop, \emph{\english{Foundations of constructive analysis}}, New York: McGraw-Hill, 1967; réimpression, New York: Ishi Press International, 2012.

E.~Bishop et D.~Bridges, \emph{\english{Constructive analysis}}, Berlin: Springer,~1985.

M.~Bridger, \emph{\english{Real analysis: a constructive approach}}, Hoboken: Wiley, 2007.

D.~Bridges et F.~Richman, \emph{\english{Varieties of constructive mathematics}}, Cambridge: Cambridge University Press,~1987.

L.~Brouwer (éd. par D.~van Dalen), \emph{\english{Brouwer's Cambridge lectures on intuitionism}}, Cambridge: Cambridge University Press, 1981.

A.~Cohen, H.~Cuypers et H.~Sterk (éds.), \emph{\english{Some tapas of computer algebra}}, New York: Springer, 1999.

D. A.~Cox, J.~Little et D.~O'Shea, \emph{\english{Ideals, varieties, and algorithms: an introduction to computational algebraic geometry and commutative algebra}}, 4\ieme\ édition, Cham: Springer, 2015.

R.~David, K.~Nour et C.~Raffalli, \emph{Introduction à la logique: théorie de la démonstration}, 2\ieme\ édition, Paris: Dunod, 2004.

G.~Dowek, \emph{La logique}, Paris: Flammarion, 1995.

S.~Feferman, \emph{\english{In the light of logic}}, Oxford: Oxford University Press, 1998.

H.~Friedman, «\english{Classically and intuitionistically provably recursive functions in Peano}», dans \emph{\english{Higher set theory: proceedings, Oberwolfach, Germany, April 13-23, 1977}} (éd.\ par G. H.~Müller et D. S.~Scott), Berlin: Springer, 1978, p.~21-27.

R.~Goodstein, \emph{\english{Recursive number theory}}, Amsterdam: North-Holland, 1957.

D.~Hilbert, «\german{Über die Grundlagen der Logik und der Arithmetik}», dans \emph{\german{Verhandlungen des dritten internationalen Mathematiker-Kongresses in Heidelberg vom 8. bis 13. August 1904}}, Leipzig: Teubner, 1905, p.~174-185, \url{http://www.mathunion.org/ICM/ICM1904}. Traduction par P.~Boutroux: «Sur les fondements de la logique et de l'arithmétique», \emph{L'enseignement mathématique}, vol.~7, 1905, p.~89-103, \doi{10.5169/seals-8424}.

D.~Knuth, «\english{Algorithmic thinking and mathematical thinking}», \emph{American Mathematical Monthly}, vol.~92, \no3, 1985, p.~170-181.

J.-L.~Krivine, \emph{Théorie des ensembles}, Paris: Cassini, 1998.

I.~Lakatos, \emph{Preuves et réfutations: essai sur la logique de la découverte mathématique}, Paris: Hermann, 1984, traduit par N.~Balacheff et J.-M.~Laborde.

R.~Lalement, \emph{Logique, réduction, résolution}, Paris: Masson, 1990.

H.~Lombardi et C.~Quitté, \emph{Algèbre commutative: méthodes constructives}, Paris: Calvage \& Mounet, 2011.

R.~Mines, F.~Richman et W.~Ruitenburg, \emph{\english{A course in constructive algebra}}, New York: Springer, 1988.

E.~Nagel, J. R.~Newman, K.~Gödel et J.-Y.~Girard, \emph{Le théorème de Gödel}, Paris: Seuil, 1989.

B.~Pire, «Hilbert (Problèmes de)», dans \emph{Encyclopædia universalis}, vol.~11, 1989, p.~412-418.

\begin{sloppypar}
  \expandafter\def\expandafter\UrlBreaks\expandafter{\UrlBreaks
  \do\e}

H.~Poincaré, \emph{Dernières pensées}, Paris: Flammarion, 1913, \url{http://henripoincarepapers.univ-lorraine.fr/chp/text/hp1913dp.html}.
\end{sloppypar}

A.~Robinson, «\english{Formalism 64}», dans \emph{\english{Logic, methodology and philosophy of science: proceedings of the 1964 International Congress}} (éd. par Y.~Bar-Hillel), Amsterdam: North-Holland, 1965, p.~228-246.

E.~Schechter, \emph{\english{Handbook of analysis and its foundations}}, New York: Academic Press, 1997.

\begin{sloppypar}
A.~Seidenberg, «\english{Constructions in algebra}», \emph{\english{Transactions of the American Mathematical Society}}, vol.~197, 1974, p.~273-313, \doi{10.1090/S0002-9947-1974-0349648-2}.
\end{sloppypar}

------, «\english{What is noetherian?}», \emph{Rendiconti del Seminario Matematico e Fisico di Milano}, vol.~44, 1974, p.~55-61.

G.~Toraldo di Francia (éd.), \emph{\italian{L'infinito nella scienza}}, Rome: \italian{Istituto della Enciclopedia Italiana}, 1987.

J.~van Heijenoort (éd.), \emph{\english{From Frege to Gödel: a source book in mathematical logic, 1879-1931}}, Cambridge: Harvard University Press, 1967.

\begin{sloppypar}
H.~Weyl, \emph{\german{Das Kontinuum: kritische Untersuchungen über die Grundlagen der Analysis}}, Leipzig: Veit \& Comp.,
1918, \url{https://archive.org/details/daskontinuumkrit00weyluoft}. Traduction par J.~Largeault dans \emph{\emph{Le continu} et autres écrits}, Paris: Vrin, 1994.
\end{sloppypar}

}\section*{Au sujet de la théorie de la mesure, une note technique.}

La théorie de la mesure a pour ambition de généraliser la théorie des grandeurs héritée des mathématiques grecques. Après les succès du calcul différentiel et intégral dans le calcul des longueurs, aires et volumes de toutes sortes de courbes, surfaces et objets solides, des problèmes de fondements sont apparus lorsqu'on a introduit, à la fin du 19\ieme\ siècle, sous l'impulsion de Cantor, Dedekind et Hilbert, de nouveaux objets beaucoup plus généraux, comme des parties arbitraires de la droite, du plan ou de l'espace.

L'ambition d'attribuer une mesure, analogue à la longueur d'un segment, à toute partie de la droite réelle, a été l'objet de tentatives multiples, dont les plus marquantes ont été la définition de la mesure des boréliens par Borel et la théorie de l'intégrale de Lebesgue.

On doit définir ce qu'est une partie mesurable, et ce qu'est sa mesure, lorsque du moins la partie n'est pas trop grande. En particulier toute partie mesurable bornée doit avoir une mesure finie.

L'idée centrale pour fonder une théorie qui tienne la route est d'imposer la «$\sigma$-additivité»: si l'on a une réunion dénombrable de parties mesurables disjointes, que chacune de ces parties a une mesure finie, et que la somme de ces mesures finies est elle-même finie, alors cette somme doit être la mesure de la réunion dénombrable. Les parties mesurables de la droite réelle doivent admettre comme cas particuliers les segments, elles doivent être stables par réunion et intersection dénombrable ainsi que par passage au complémentaire. Borel définit ce que nous appelons aujourd'hui les boréliens comme les parties de la droite qui peuvent être construites selon les processus invoqués ci-dessus à partir des segments ouverts. Et il démontre que l'on peut attribuer une mesure finie à tout borélien contenu dans un segment fermé borné en respectant la propriété de $\sigma$-additivité. Cette magnifique théorie a été reprise par Bishop de manière entièrement constructive, mais elle a été éclipsée par l'intégrale de Lebesgue. C'est la théorie de Lebesgue qui est aujourd'hui couramment enseignée à l'université, en troisième année de licence ou première année de master. La théorie de Borel est plutôt considérée comme une théorie de spécialistes, trop compliquée pour être enseignée dans le cursus usuel.

De manière générale, toutes les questions concernant la théorie de la mesure présentent des défis intéressants.

Un\textperiodcentered e mathématicien\textperiodcentered ne classique donne souvent l'argument suivant: si l'on ne considère que les nombres réels récursifs, ils forment un ensemble dénombrable, donc de mesure nulle, et les mathématiques usuelles s'effondrent.  Vouloir se restreindre aux êtres mathématiques explicites est donc une catastrophe.

Encore un argument qui oublie la différence entre Vercingétorix et Don Quichotte!

En fait, en théorie de la mesure classique, une partie mesurable de $[0,1]$ de mesure strictement positive peut très bien ne contenir aucun réel que l'on soit capable d'expliciter. Ceci devrait plutôt être considéré comme une faiblesse de la théorie classique que comme une force. Et c'est un résultat qu'en général on cache soigneusement aux étudiant\textperiodcentered e\textperiodcentered s pour ne pas les décourager. On leur demande d'admettre que l'ensemble $\mathbb{R}$ peut être muni d'une relation de bon ordre, sans trop insister sur l'étrangeté de la chose, mais on ne leur parle guère du «paradoxe» apparent que nous venons d'indiquer en théorie de la mesure.

Quand on lit Bishop au contraire, on voit que non seulement il est impossible d'énumérer les réels, conformément à Cantor, mais aussi que toute partie mesurable de mesure strictement positive contient vraiment des réels explicites\footnote{En fait on peut expliciter un ensemble de Cantor dans toute partie mesurable de mesure strictement positive. Une preuve constructive dans un cas particulier est donnée dans le livre \emph{Épistémologie mathématique} (H.~Lombardi, Paris: Ellipses, 2011), section~8.3.} que l'on sait calculer. Ceci est naturellement un point très positif en faveur de l'approche minimale de Bishop. Et il n'y a là aucune contradiction flagrante avec le «paradoxe» classique précédent. En effet, force recours au tiers exclu est nécessaire aux mathématicien\textperiodcentered ne\textperiodcentered s classiques pour énumérer les réels récursifs. Du point de vue de Bishop, il est tout simplement impossible de montrer que les réels récursifs forment une partie mesurable dans l'intervalle $[0,1]$.  À fortiori ils ne forment pas une partie que l'on pourrait énumérer de manière effective.

Maintenant, la preuve en mathématiques minimales, par Bishop, que tout borélien explicite de mesure strictement positive contient une infinité non dénombrable de points explicites est aussi une preuve en mathématiques classiques que tout borélien de mesure strictement positive contient une infinité non dénombrable de points. Il suffit en effet de relire la preuve de Bishop en vidant de son sens intuitif le mot «explicite» et en décrétant en son for intérieur que, vu le principe du tiers exclu, on est autorisé à admettre que tout est \emph{explicite} en mathématiques, mais en un sens affaibli.

En fait, à quoi sert vraiment la théorie de la mesure? à mettre en valeur des bizarreries purement formelles sans réalité concrète comme le paradoxe précédent\footnote{Ou encore le fameux paradoxe de Banach-Tarski?}?  ou à obtenir des résultats concrets qui ont une signification claire pour tout le monde?

Vers la fin de l'article, Apéry rappelle le problème ouvert suivant.
\begin{quotation}
\noindent Pour presque tout réel~$\alpha>1$, c'est-à-dire sauf sur un ensemble de mesure nulle, les~$\alpha^n$ sont «bien répartis» sur le groupe additif de~$\mathbb{R}$/$\mathbb{Z}$; néanmoins, un problème important et non résolu est de nommer un~$\alpha$ tel que les~$\alpha^n$ soient bien répartis.
\end{quotation}

D'après les spécialistes du sujet, le problème semble avoir été résolu positivement en russe dès 1961 et en anglais vers 1996. Un livre de référence sur les questions de distribution de réels modulo~1 est le suivant.\medskip

{\leftskip15pt
\parindent-15pt

Y.~Bugeaud, \emph{Distribution modulo one and Diophantine approximation}, New York: Cambridge University Press, 2012.\medskip

}L'auteur nous a recommandé de lire la page~46. Citons l'extrait suivant (nous traduisons en français).
\begin{quotation}
Levin [430] (voir aussi [429, 431] et Kulikova [412]) ont construit des nombres réels~$\alpha$ tels que la suite des~$\alpha^n$ est distribuée uniformément modulo~1 [\dots]. Dans un manuscrit non publié, Lerma [424] a donné une construction alternative très compliquée de nombres réels ayant cette même propriété.
\end{quotation}

Les références citées par Bugeaud sont les suivantes.\medskip

{\setbox0\hbox{[412]\enskip}\leftskip\wd0
\parindent-\wd0

[412]\enskip M. F.~Kulikova, «\english{A construction problem concerned with the distribution of the fractional parts of an exponential function}», \emph{Doklady Akademii Nauk SSSR}, vol.~143, 1962, p.~522-524 (en russe). Traduction anglaise: \emph{\english{Soviet Mathematics: Doklady}}, vol.~3, 1962, p.~422-425.

\begin{sloppypar}
[424]\enskip M. A.~Lerma, «\english{Construction of a number greater than one whose powers are uniformly distributed modulo one}», \emph{prépublication}, 1996, \url{http://sites.math.northwestern.edu/~mlerma/papers/constr_ud_mod1.pdf}.
\end{sloppypar}

[429]\enskip M.~Levin, «\english{Completely uniform distribution of fractional parts of the exponential function}», \emph{Trudy Seminara im. I. G.~Petrovskogo}, vol.~7, 1981, p.~245-256 (en russe). Traduction anglaise: \emph{Journal of Soviet Mathematics}, vol.~31, 1985, p.~3247-3256.

[430]\enskip ------, «\textcyrillic{\CYREREV\cyrf\cyrf\cyre\cyrk\cyrt\cyri\cyrv\cyri\cyrz\cyra\cyrc\cyri\cyrya\ \cyrt\cyre\cyro\cyrr\cyre\cyrm\cyrery\ \CYRK\cyro\cyrk\cyrs\cyrm\cyra} 
[Effectivisation du théorème de Koksma]», \emph{Matematicheskie Zametki}, vol.~47, 1990, p.~163-166.

\begin{sloppypar}
[431]\enskip ------, «\english{Jointly absolutely normal numbers}», \emph{Matematicheskie Zametki}, vol.~48, 1990, p.~61-71 (en russe). Traduction anglaise: \emph{Mathematical Notes}, vol.~48, 1990, p.~1213-1220.\medskip
\end{sloppypar}

}Il est dommage que la référence [430] n'ait pas fait l'objet d'une traduction. Nous ne savons pas si les démonstrations fournies sont vraiment constructives, mais il est probable qu'elles le soient.
\[*\ *\ *\]
\end{document}